\documentclass[reqno,11pt]{article}

\usepackage{amssymb,amsmath}
\usepackage{enumitem}

\usepackage{epsfig}
\usepackage[normal]{subfigure}
\usepackage{graphicx,color}

\usepackage[latin1]{inputenc}

\newtheorem{theo}{Theorem}
 \newtheorem{lem}[theo]{Lemma}
\newtheorem{prop}[theo]{Proposition}
\newtheorem{defn}[theo]{Definition}
\newtheorem{rem}{Remark}

\def\lip{{\rm Lip}}

\def\E{\mathbb{E}}

\def\C{\mathbb{C}}

\def\phi{\varphi}

\newcommand {\nn}{\nonumber}

\newcommand {\noi}{\noindent}

\def\mathsf{\bf}
\def\N{\mathbb{N}}

\def\R{\mathbb{R}}
\def\Z{\mathbb{Z}}

\def\E{\mathrm E}

\def\text{\mbox}

\def\1{{\bf 1}}

\newcommand\beqn{\begin{displaymath}}  
\newcommand\eeqn{\end{displaymath}}

\begin{document}

\title{A generalized nonlinear model for long memory
conditional heteroscedasticity }
\author{Ieva Grublyt\.e$^{1,2}$  and  Andrius \v Skarnulis$^{2}$  }
\date{\today \\  \small
\vskip.2cm
$^1$ Universit\'e de Cergy Pontoise, D\'epartament de Math\'ematiques, 95302 Cedex, France \\
$^2$ Vilnius University, Institute of Mathematics and Informatics, 08663 Vilnius, Lithuania}
\maketitle

\begin{abstract} We study the existence and properties of stationary solution of ARCH-type equation
$r_t= \zeta_t \sigma_t$, where $\zeta_t$ are standardized i.i.d. r.v.'s and the conditional
variance satisfies an AR(1) equation $\sigma^2_t = Q^2\big(a + \sum_{j=1}^\infty b_j  r_{t-j}\big) + \gamma \sigma^2_{t-1}$ with a
Lipschitz function $Q(x)$ and real parameters $a, \gamma, b_j $. The paper extends the model and the results in
\cite{dou2015} from the case $\gamma =  0$ to the case $0< \gamma < 1$. We also obtain a new  
condition for the existence of higher moments of $r_t$ which does not include the Rosenthal 
constant. 
In the particular case when
$Q$ is the square root of a quadratic polynomial, we prove that $r_t$ can exhibit a leverage effect and long memory. We also
present simulated trajectories and histograms of marginal density of $\sigma_t$ for different values of
$\gamma$.

\end{abstract}

\begin{quote}


\end{quote}

\section{Introduction}

Doukhan et al.~\cite{dou2015} discussed the existence of stationary solution of
conditionally heteroscedastic equation
\begin{eqnarray}\label{genform0}
r_t&=&\zeta_t \sigma_t, \qquad
\sigma^2_t\ =\  Q^2\Big(a + \sum_{j=1}^\infty b_j  r_{t-j}\Big),
\end{eqnarray}
where $\{\zeta_t\}$ are standardized i.i.d. r.v.'s, $a, b_j $ are real parameters and  $Q(x)$ is a Lipschitz function of real variable $x \in \R$. Probably, the most important
case of (\ref{genform0}) is
\begin{equation}\label{Qform1}
Q(x) = \sqrt{c^2 +  x^2},
\end{equation}
where $c \ge 0$ is a parameter.  The model (\ref{genform0})-\eqref{Qform1} includes the classical Asymmetric
ARCH(1) of Engle \cite{eng1990} and the Linear ARCH (LARCH) model  of Robinson \cite{rob1991}:
\begin{equation}\label{larch}
r_t\ =\  \zeta_t \sigma_t, \qquad
\sigma_t\ =\  a + \beta \sum_{j=1}^\infty j^{d-1}  r_{t-j}.
\end{equation}
\cite{gir2000} proved that the squared stationary solution  $\{r^2_t\}$ of the LARCH model in \eqref{larch}
with $b_j$ decaying as $j^{d-1}, 0< d< 1/2 $
may have long memory autocorrelations.
The leverage effect in the LARCH
model was discussed in detail in \cite{gir2004}. Doukhan et al. \cite{dou2015} extended the above properties
of the LARCH model (long memory and leverage) to
the model in (\ref{genform0})-\eqref{Qform1} with $c>0$ or strictly positive volatility.

The present paper extends the results of \cite{dou2015} to a more general class of volatility forms:
\begin{eqnarray}\label{genform}
r_t&=&\zeta_t \sigma_t, \qquad
\sigma^2_t\ =\ Q^2\Big(a + \sum_{j=1}^\infty b_j  r_{t-j}\Big) + \gamma \sigma^2_{t-1},
\end{eqnarray}
where $\{\zeta_t\}, a, b_j, Q(x) $ are as in (\ref{genform0}) and $0 < \gamma < 1 $ is a parameter.  The inclusion
of lagged $\sigma^2_{t-1}$ in (\ref{genform}) helps to reduce very sharp peaks and  clustering of volatility which
occur in trajectory of (\ref{genform0})-\eqref{Qform1}
near the threshhold $c>0$ (see Fig. 1).  The generalization from (\ref{genform0}) to (\ref{genform})
is similar to that from ARCH to GARCH models, see \cite{eng1982}, \cite{bol1986}, particularly,  (\ref{genform}) with $Q(x) $ of $\eqref{Qform1}$ and
$b_j = 0, j\ge 2 $ reduces
to the Asymmetric GARCH(1,1) of  Engle \cite{eng1990}.

Let us describe the main results of this paper. 
Sec.~2 (Theorems~\ref{Xexists}  and \ref{Xeven}) obtain sufficient conditions for the existence
of stationary solution of (\ref{genform}) with $\E |r_t|^p < \infty $ and
$\gamma \in [0,1)$. Theorem~\ref{Xexists} extends the corresponding result 
in (\cite{dou2015}, Thm.~4) from $\gamma =0$ to $\gamma > 0$. 
Theorem~\ref{Xeven} is new even in the case $\gamma =0$ by providing an explicit sufficient condition \eqref{Mp}
for higher-order even moments ($p=4,6,\dots$) which does not involve the absolute constant in the Burkholder-Rosenthal inequality \eqref{rosen}.
Condition~\eqref{Mp} coincides with the corresponding moment condition for the LARCH model and
is important for statistical applications, 
see Remark \ref{remMp}. 
The remaining sec.~3-5 deal exclusively with the case of quadratic
$Q^2$ in \eqref{Qform1}, referred to as the Generalized Quadratic ARCH (GQARCH) model in the sequel. 
Theorem \ref{long} (sec.~3) obtains long memory properties of the squared process $\{r^2_t\} $ of
the GQARCH model with $\gamma \in (0,1)$ and coefficients
$b_j$ decaying regularly as $b_j \sim \beta j^{d-1}, \, j \to \infty, \, 0< d < 1/2$. Similar properties
were established in \cite{dou2015} for the GQARCH model with $\gamma =0$ and for the LARCH model
\eqref{larch} in \cite{gir2000}, \cite{gir2004}. Quasi-maximum likelihood estimation for parametric GQARCH model with long memory was recently studied in \cite{gru2015}. 
See the review paper \cite{gir2009} and recent work \cite{gir2014}
on long memory ARCH modeling. 
Sec.~4 extends to the GQARCH model
the leverage effect discussed in \cite{dou2015} and \cite{gir2004}. Sec.~5 presents some
simulations and volatility profiles for the LARCH and GQARCH models with parameters estimated from real data.
A general impression from our results
is that the GQARCH modification in (\ref{genform}), \eqref{Qform1}
of the QARCH model in \cite{dou2015}
allows for a more realistic volatility modeling as compared to the LARCH and QARCH  models,
at the same time preserving
the long memory and the leverage properties of the above mentioned models.

\section{Stationary solution}
Denote $|\mu|_p := \E |\zeta_0|^p \,  (p>0), \, \mu_p := \E \zeta_0^p \, (p=1,2, \dots)$ and let 
\begin{equation}\label{Xform}
X_t \ := \ \sum_{s< t} b_{t-s} r_s.
\end{equation}
Since $0\le \gamma  < 1 $, equations \eqref{genform} yield
\begin{eqnarray}
\sigma^2_t&=&\sum_{\ell =0}^\infty \gamma^\ell Q^2(a + X_{t-\ell}) \quad \text{and} \quad
\label{genformQ}
r_t\ =\  \zeta_t \sqrt{ \sum_{\ell =0}^\infty \gamma^\ell Q^2(a + X_{t-\ell})}.
\end{eqnarray}
In other words, stationary solution of   (\ref{genform}), or
\begin{equation} \label{rformQ}
r_t \ = \  \zeta_t \sqrt{ \sum_{\ell =0}^\infty \gamma^\ell Q^2(a + \sum_{j=1}^\infty b_j r_{t-\ell -j} )}
\end{equation}
can be defined via \eqref{Xform}, or stationary solution of
\begin{equation}\label{XformQ}
X_t \ := \ \sum_{s< t} b_{t-s} \zeta_s \sqrt{ \sum_{\ell =0}^\infty \gamma^\ell Q^2(a + X_{s-\ell})},
\end{equation}
and vice versa.

\noindent
In Theorem \ref{Xexists}  below, we assume that $Q$ in \eqref{genformQ} is a
Lipschitz function, i.e., there exists $ \lip_Q >0$ such that
\begin{equation}\label{QLip}
|Q(x)-Q(y)| \le \lip_Q |x-y|, \qquad  x,y \in \R.
\end{equation}
Note \eqref{QLip} implies the bound
\begin{equation}\label{Qnel}
Q^2(x) \le c_1^2 + c_2^2 x^2, \qquad x \in \R,
\end{equation}
where $c_1\ge 0, \, c_2 \ge  \lip_Q$ and $c_2 $ can be chosen arbitrarily close to $\lip_Q$.

Let us give some formal  definitions. Let ${\cal F}_t = \sigma(\zeta_s, s\le t), t \in \Z$ be the sigma-field generated
by $\zeta_s, s\le t$. A random process $\{ u_t, t \in \Z\}$ is called {\it adapted} (respectively,  {\it predictable})
if $u_t$ is ${\cal F}_t$-measurable for each $t \in \Z $ (respectively,  $u_t$ is ${\cal F}_{t-1}$-measurable for each $t \in \Z $).

\begin{defn} Let $p>0$ be arbitrary.

\smallskip

\noi (i) By $L^p$-solution of   \eqref{genformQ} or/and \eqref{rformQ}
we mean  an adapted process $\{r_t, t \in \Z\}$
with $\E |r_t|^p < \infty $ such that  for any $t \in \Z$ the series
$X_t = \sum_{j=1}^\infty  b_j r_{t-j}$ converges in $L^p$, the series
$\sigma^2_t = \sum_{\ell=0}^\infty \gamma^\ell Q^2(a+ X_{t-\ell}) $
converges in $L^{p/2}$ and   \eqref{rformQ} holds.

\smallskip

\noi (ii) By $L^p$-solution of \eqref{XformQ} we mean  a predictable process $\{X_t, t \in \Z\}$ with
$\E |X_t|^p < \infty $ such that  for any $t \in \Z$ the series $\sigma^2_t = \sum_{\ell=0}^\infty \gamma^\ell Q^2(a+ X_{t-\ell})$
converges in $L^{p/2}$, the series
$\sum_{s<t}  b_{t-s} \zeta_s \sigma_s$ converges in $L^p$ and \eqref{XformQ}  holds.
\end{defn}

Define
\begin{equation*}\label{Bpdef}
B_p := \begin{cases}
\sum_{j=1}^\infty |b_j|^p, &0<p < 2, \\
\big(\sum_{j=1}^\infty b_j^2\big)^{p/2}, &p\ge 2,
\end{cases}
\qquad
B_{p,\gamma} := \begin{cases}
B_p/(1-\gamma^{p/2}), &0<p < 2, \\
B_p/(1-\gamma)^{p/2}, &p\ge 2.
\end{cases}
\end{equation*}
Note $B_p = B_{p,0}$.
As in \cite{dou2015}, we use the following moment inequality, see \cite{bur1973}, \cite{von1965},  \cite{ros1970}.

\begin{prop} \label{Yp} Let
$\{Y_j, j \ge 1\}$ be a sequence of r.v.'s such that $\E |Y_j|^p < \infty $ for some
$p >0$. If $p>1 $ we additionally assume that $\{Y_j\}$ is a martingale difference sequence:
$\E [Y_j |Y_1, \dots, Y_{j-1}] = 0, \, j=2,3, \dots $. Then
there exists a constant $K_p\ge 1$ depending only on $p$ and such that
\begin{equation}\label{rosen}
\E \big|\sum_{j=1}^\infty Y_j\big|^p \ \le \ K_p \begin{cases}
\sum_{j=1}^\infty \E |Y_j|^p, &0< p \le 2, \\
\big(\sum_{j=1}^\infty (\E |Y_j|^p)^{2/p}\big)^{p/2}, &p > 2.
\end{cases}
\end{equation}

\end{prop}

\smallskip

Proposition \ref{Xreq} says that equations \eqref{rformQ} and \eqref{XformQ} are equivalent in the sense that
by solving one the these  equations one readily obtains a solution to the other one.

\begin{prop} \label{Xreq} Let $Q$ be a measurable function satisfying  \eqref{Qnel} with some $c_1, c_2 \ge 0$ and $\{\zeta_t\}$ be an i.i.d. sequence with $|\mu|_p = \E|\zeta_0|^p < \infty$ and satisfying
$\E \zeta_0 = 0$ for $p> 1$. In addition, assume $B_p < \infty$ and $0\le \gamma < 1 $.

\smallskip

\noi (i) Let $\{X_t \}$ be a stationary $L^p$-solution of \eqref{XformQ} and let
$\sigma_t := \sqrt{ \sum_{\ell =0}^\infty \gamma^\ell Q^2(a + X_{t-\ell})}$.   Then $\{r_t =
\zeta_t \sigma_t  \}$ in \eqref{genformQ}  is a stationary
$L^p$-solution of \eqref{rformQ} and
\begin{equation}\label{rX}
\E |r_t|^p \ \le \ C(1+ \E |X_t|^p).
\end{equation}
Moreover, for $p > 1 $, $\{r_t, {\cal F}_t, t \in \Z\}$ is a
martingale difference sequence with
\begin{equation}\label{genformvol}
\E [r_t|{\cal F}_{t-1}] = 0, \qquad \E [|r_t|^p| {\cal F}_{t-1}]  = |\mu|_p \sigma_t^p.
\end{equation}

\noi (ii) Let $\{r_t \}$ be a stationary $L^p$-solution of \eqref{rformQ}.
Then  $\{X_t\}$ in \eqref{Xform} is a   stationary $L^p$-solution of \eqref{XformQ} such that
\begin{equation*}\label{Xr}
\E |X_t|^p \ \le \ C\E |r_t|^p.
\end{equation*}
Moreover, for $p \ge 2$
\begin{equation*}\label{covX}
\E [X_t X_0] \ = \ \E r^2_0 \sum_{s=1}^\infty b_{t+s} b_s, \qquad t = 0,1,\dots.
\end{equation*}

\end{prop}

\noi {\it Proof.} (i)
First, let $0< p \le 2 $. Then
$\E |\sigma_t|^p =  \E |\sigma^2_t|^{p/2} \le   \sum_{\ell=0}^\infty |\gamma^{p/2}|^\ell \E |Q(a+ X_{t-\ell})|^p < \infty. $
Hence, using \eqref{Qnel}, the fact that $\{X_t\}$ is predictable
and $|Q(a + X_{t-\ell})|^p \le
|c_1^2 + c_2^2(a+ X_{t-\ell})|^{p/2} \le C(1+ |a+ X_{t-\ell}|^p) \le C(1 + |X_{t-\ell}|^p) $
we obtain
\begin{eqnarray*}
\E |r_t|^p
&=&|\mu|_p\E |\sigma_t|^p \  \le \  C \sum_{\ell=0}^\infty |\gamma^{p/2}|^\ell (1 + \E |X_{t-\ell}|^p)\\
&\le&C(1+ \E |X_t|^p) \ < \ \infty,
\end{eqnarray*}
proving \eqref{rX} for $p \le 2$.  Next, let $p> 2 $. Then $\E |\sigma_t|^p \le (\sum_{\ell=0}^\infty \gamma^\ell \E^{2/p} |Q(a+X_t)|^p )^{p/2} \le
C \E |Q(a+ X_t)|^p$
by stationarity and Minkowski's inequality and hence \eqref{rX} follows using the same argument as above. Clearly, for $p  > 1 $
$\{r_t = \zeta_t \sigma_t\}$ is a martingale difference sequence and satisfies \eqref{genformvol}.
Then, the convergence in $L^p$ of
the series in \eqref{Xform} follows from \eqref{rX} and Proposition  \ref{Yp}:
\begin{eqnarray*}
\E |\sum_{j=1}^\infty b_j r_{t-j}|^p &\le&C\left\{\begin{array}{ll}
\sum_{j=1}^\infty |b_j|^p, &0< p \le 2 \\
\big(\sum_{j=1}^\infty b_j^2\big)^{p/2}, &p > 2
\end{array}
\right\} \ = \ CB_p \ < \ \infty.
\end{eqnarray*}
In particular,
$\zeta_t \sqrt{\sum_{\ell =0}^\infty \gamma^\ell Q^2(a + \sum_{s<t} b_{t-\ell- s} r_s)} = \zeta_t \sqrt{\sum_{\ell=0}^\infty \gamma^\ell
Q^2(a+ X_{t-\ell})} = r_t $
by the definition of $r_t$.  Hence, $\{r_t \}$ is a $L^p$-solution of (\ref{rformQ}).  Stationarity of   $\{r_t \}$ follows
from stationarity of $\{X_t \}$.

\smallskip

\noi (ii) Since $\{r_t\}$ is a $L^p$-solution  of \eqref{rformQ}, so $r_t = \zeta_t \sigma_t = \zeta_t \sqrt{\sum_{\ell=0}^\infty
\gamma^\ell
 Q^2(a + X_{t-\ell})}$ with $X_t$ defined
in \eqref{Xform} and $\{X_t\}$ satisfy \eqref{Xform}, where the series converges in $L^p$. The rest follows
as in \cite{dou2015}, proof of Prop.3.   \hfill $\Box$

\begin{rem} \label{Lpsol} {\rm  Let $p\ge 2$ and $|\mu|_p < \infty$, then  by inequality \eqref{rosen},
$\{r_t\}$ being a stationary $L^p$-solution of \eqref{genformQ} is equivalent
to  $\{r_t\}$ being a stationary $L^2$-solution of \eqref{genformQ} with $\E |r_0|^p < \infty $.
Similarly, if $Q$ and $\{\zeta_t\}$  satisfy the conditions of Proposition \ref{Xreq} and $p\ge 2$,
then  $\{X_t\}$ being a stationary $L^p$-solution of \eqref{Xform} is equivalent
to  $\{X_t\}$ being a stationary $L^2$-solution of \eqref{Xform} with $\E |X_0|^p < \infty $. 
See also (\cite{dou2015}, Remark 1). 
}
\end{rem}

\begin{theo} \label{Xexists} Let $\{\zeta_t\}$ satisfy the conditions of Proposition \ref{Xreq} and
$Q$ satisfy the Lipschitz condition in \eqref{QLip}.

\medskip

\noi (i) Let $p>0$ and
\begin{equation} \label{cQB}
K^{1/p}_p\, |\mu|^{1/p}_p \,\lip_Q\, B^{1/p}_{p,\gamma}   < 1,
\end{equation}
where $K_p$ is the absolute constant from the moment inequality in \eqref{rosen}.
Then there exists a unique stationary $L^p$-solution $\{X_t\}$ of \eqref{XformQ} and
\begin{equation}\label{X2mom}
\E |X_t|^p \ \le \frac{C(p,Q) |\mu|_p  B_p }{1 - K_p |\mu|_p \lip_Q^p B_{p,\gamma}},
\end{equation}
where $C(p,Q) < \infty $ depends only on $p$ and $c_1, c_2 $ in \eqref{Qnel}.

\smallskip

\noi (ii) Assume, in addition, that
$Q^2(x) = c_1^2 + c_2^2 x^2$, where $c_i \ge 0, i=1,2$, and $\mu_2 = \E \zeta^2_0 =1 $.
Then $c^2_2 B_{2,\gamma} < 1$ is a necessary and sufficient condition for the existence
of  a stationary $L^2$-solution $\{X_t\}$ of \eqref{XformQ} with $a \neq 0$.
\end{theo}

\noi {\it Proof.} (i) We follow the proof of Theorem~4 in \cite{dou2015}. For $n\in \N $ we recurrently define a solution 
of  \eqref{XformQ} with zero initial condition at $t \le - n $ as
\begin{eqnarray}
X^{(n)}_t&:=&\begin{cases}
0, &t \le -n, \\
\sum_{s= -n}^{t-1}b_{t-s} \zeta_s \sigma^{(n)}_s, &t> -n, \quad t \in \Z,
\end{cases}\label{Xndef}
\end{eqnarray}
where $ \sigma^{(n)}_s :=  \sqrt{\sum_{\ell =0}^{n+s} \gamma^\ell Q^2(a + X^{(n)}_{s-\ell})}$.
Let us show that $\{X^{(n)}_t \}$ converges in $L^p$ to a stationary $L^p$-solution $\{X_t\}$
as $n \to \infty $.

\noi First, let $0< p \le 2$.
Let $m > n \ge 0$.  Then by inequality \eqref{rosen} for any $ t>  -m$ we have that
\begin{eqnarray}
\E |X^{(m)}_t - X^{(n)}_t|^p&\le&K_p|\mu|_p \Big\{ \sum_{-m \le s<-n}  |b_{t-s}|^p \E |\sigma^{(m)}_s|^p
+ \sum_{-n \le s< t}  |b_{t-s}|^p \E |\sigma^{(m)}_s - \sigma^{(n)}_s|^p \Big\} \nn \\
&&=:K_p|\mu|_p \big\{S'_{m,n} + S''_{m,n} \big\}.  \label{sumS}
\end{eqnarray}
Using $|Q(a+x)|^p \le C + c_3^p |x|^p, x \in \R $ with $c_3 > c_2 > \lip_Q$ arbitrarily close to $\lip_Q$, see \cite{dou2015}, proof of Theorem~4,
we obtain
\begin{eqnarray}\label{S1mn}
S'_{m,n}&\le&\sum_{-m \le s<-n}  |b_{t-s}|^p \sum_{\ell =0}^{m+s}\gamma^{p\ell/2} (C + c_3^p \E |X^{(m)}_{s-\ell}|^p) .
\end{eqnarray}
Next, using $| (\sum_{i>0} \gamma^i x_i^2)^{1/2} -  (\sum_{i>0} \gamma^i y_i^2)^{1/2}| \le
(\sum_{i>0} \gamma^i (x_i-y_i)^2)^{1/2}$ we obtain
\begin{equation}\label{sigmabdd}
|\sigma^{(m)}_s - \sigma^{(n)}_s| \le \left(\sum_{\ell=0}^{s+n} \gamma^\ell
\big(Q(a + X^{(m)}_{s-\ell}) -Q(a + X^{(n)}_{s-\ell}) \big)^2+
\sum_{\ell=s+n+1}^{s+m} \gamma^\ell Q^2(a + X^{(m)}_{s-\ell}) \right)^{1/2}.
\end{equation}
Hence from the Lipschitz condition in \eqref{QLip} we have that
\begin{eqnarray*}
S''_{m,n}
&\le&\sum_{-n \le s<t}  |b_{t-s}|^p
\left(\sum_{\ell=0}^{s+n} \gamma^{p \ell /2}
\lip_Q^p  \E | X^{(m)}_{s-\ell}- X^{(n)}_{s-\ell} |^p+
\sum_{\ell=s+n+1}^{s+m} \gamma^{p \ell /2}(C+ c_3^p\E | X^{(m)}_{s-\ell})  |^{p})\right).
\end{eqnarray*}
Combining \eqref{sumS}  and the above bounds we obtain
\begin{eqnarray}
\E |X^{(m)}_t - X^{(n)}_t|^p&\le&
K_p |\mu|_p \Big(
c_3^p \sum_{-m \le s<t}  |b_{t-s}|^p \sum_{\ell=0}^{s+m} \gamma^{p \ell /2}
  \E | X^{(m)}_{s-\ell}- X^{(n)}_{s-\ell} |^p\nn  \\
  &+&\quad C\sum_{-m \le s<-n}  |b_{t-s}|^p
\sum_{\ell=0}^{s+m} \gamma^{p \ell /2}
+ C\sum_{-n \le s<t}  |b_{t-s}|^p
\sum_{\ell=s+n+1}^{s+m} \gamma^{p \ell /2}\Big) \nn \\
&\le&C K_p|\mu|_p \kappa^p_{t+n,\gamma}
+  K_p|\mu|_p c_3^p  \sum_{-m \le s<t}  b_{t-s,\gamma}^p  \E |X^{(n)}_s - X^{(m)}_s |^p.
\label{Xmn}
\end{eqnarray}
where $b_{s,\gamma}^p := \sum_{j = 0}^{s-1} \gamma^{ jp/2} |b_{s-j}|^p, s \ge 0,
\kappa_{t+n} := C (1-\gamma^{p/2
})^{-1} K_p|\mu|_p(\sum_{j>t+n} |b_j|^p + b_{t+n, \gamma}^p
) \to 0 \, (n \to \infty)$.

\smallskip

\noi Iterating inequality \eqref{Xmn} as in \cite{dou2015}, (6.55) and  using
$ K_p|\mu|_p c_3^p  \sum_{s<t}  b_{t-s,\gamma}^p =  K_p|\mu|_p c_3^p B_{p,\gamma} < 1 $ we obtain
$\lim_{m,n \to \infty} \E |X^{(m)}_t - X^{(n)}_t|^p = 0$ and hence
the existence of  $X_t$ such that  $\lim_{n \to \infty} \E |X^{(n)}_t - X_t|^p = 0$ and
satisfying the bound in \eqref{X2mom}.

\smallskip

Next, consider the case  $p> 2 $. Let $m > n \ge 0$.  Then by inequality \eqref{rosen} for any $ t>  -m$ we have that
\begin{eqnarray}
\E |X^{(m)}_t - X^{(n)}_t|^p&\le&K_p|\mu|_p \Big( \sum_{-m \le s<-n}  b^2_{t-s} \E^{2/p} |\sigma^{(m)}_s|^p
+ \sum_{-n \le s< t}  b^2_{t-s} \E^{2/p} |\sigma^{(m)}_s - \sigma^{(n)}_s|^p \Big)^{p/2} \nn \\
&&=:K_p|\mu|_p \big(R'_{m,n} + R''_{m,n} \big)^{p/2}.  \label{sumS1}
\end{eqnarray}
Similarly to \eqref{S1mn},
\begin{eqnarray*}
R'_{m,n}&\le&\sum_{s=-m}^{-n-1} b^2_{t-s}  \sum_{\ell =0}^{m+s}\gamma^\ell \E^{2/p} |Q(a +X^{(m)}_{s-\ell})|^p  \
\le\ \sum_{s=-m}^{-n-1} b^2_{t-s}  \sum_{\ell =0}^{m+s}\gamma^\ell (C + c^2_3 \E^{2/p} |X^{(m)}_{s-\ell}|^p)
\end{eqnarray*}
and, using  \eqref{sigmabdd},
\begin{eqnarray*}
R''_{m,n}
&\le&\sum_{-n \le s<t} b^2_{t-s} \E^{2/p} \Big| \sum_{\ell=0}^{s+n} \gamma^\ell (Q(a+ X^{(m)}_{s-\ell}) - Q(a+ X^{(n)}_{s-\ell}))^2
+ \sum_{\ell =s+n+1}^{s+m} \gamma^\ell Q^2(a+ X^{(m)}_{s-\ell}) \Big|^{p/2} \\
&\le&\sum_{-n \le s<t} b^2_{t-s} \Big( \sum_{\ell=0}^{s+n} \gamma^\ell \E^{2/p} |Q(a+ X^{(m)}_{s-\ell}) - Q(a+ X^{(n)}_{s-\ell})|^p
+ \sum_{\ell =s+n+1}^{s+m} \gamma^\ell \E^{2/p} |Q(a+ X^{(m)}_{s-\ell})|^p  \Big) \\
&\le&\sum_{-n \le s<t} b^2_{t-s} \Big( \lip^2_Q \sum_{\ell=0}^{s+n} \gamma^\ell \E^{2/p} |X^{(m)}_{s-\ell} - X^{(n)}_{s-\ell}|^p
+ \sum_{\ell =s+n+1}^{s+m} \gamma^\ell (C + c_3^2  \E^{2/p} |X^{(m)}_{s-\ell}|^p)  \Big)
\end{eqnarray*}
Consequently,
\begin{eqnarray*}
\E^{2/p} |X^{(m)}_t - X^{(n)}_t|^p&\le&
\kappa_{t+n} +  K^{2/p}_p|\mu|_p^{2/p} c^2_3 \sum_{-m \le s<t} b^2_{t-s,\gamma}    \E^{2/p} |X^{(m)}_{s}) - X^{(n)}_{s}|^p,
\end{eqnarray*}
where $\kappa_{t+n} := C (1-\gamma)^{-1} K^{2/p}_p|\mu|_p^{2/p} (\sum_{j>t+n} b^2_j + b_{t+n, \gamma}^2
) \to 0 \, (n \to \infty)$. 
By iterating
the  last displayed equation and using $K^{2/p}_p|\mu|_p^{2/p} c^2_3 \sum_{j=1} b^2_{j,\gamma} = K^{2/p}_p|\mu|_p^{2/p} c^2_3 B_2/(1-\gamma) < 1 $
we obtain $\lim_{m,n \to \infty} \E^{2/p} |X^{(m)}_t - X^{(n)}_t|^p = 0$ and hence
the existence of  $X_t$ such that  $\lim_{n \to \infty} \E |X^{(n)}_t - X_t|^p = 0$ and
satisfying the bound in \eqref{X2mom}. The rest of the proof of part (i)
is similar as in  \cite{dou2015}, proof of Theorem~4, and we omit the details.

\smallskip

\noi (ii) Note that $Q(x) =  \sqrt{c_1^2 + c_2^2 x^2}$ is a Lipschitz function and satisfies \eqref{QLip} with $\lip_Q = c_2$. Hence by $K_2 =1 $ and part (i), a unique
$L^2$-solution $\{X_t\}$ of \eqref{XformQ} under the condition $c^2_2 B_{2,\gamma} = c_2^2 B_2/(1-\gamma)  < 1$ exists. To show
the necessity of the last condition,
let $\{X_t \} $ be a stationary $L^2$-solution of  \eqref{XformQ}.  Then
\begin{eqnarray*}
\E X^2_t&=&\sum_{s<t}  b^2_{t-s} \sum_{\ell =0}^\infty \gamma^\ell\E Q^2(a + X_{s-\ell})\\
&=&\sum_{s<t}  b^2_{t-s} \sum_{\ell =0}^\infty \gamma^\ell \E \big(c_1^2 + c_2^2 (a+ X_{s-\ell}^2 \big) \\
&=&(B_2/(1-\gamma))\big(c_1^2 + c^2_2 (a^2 + \E X^2_t)\big) \ > c^2_2 (B_2/(1-\gamma))\E X^2_t
\end{eqnarray*}
since $a \neq 0$.  Hence, $c^2_2 B_2/(1-\gamma)  < 1$ unless $\E X^2_t = 0 $,
or  $\{ X_t = 0\}$ is a trivial process.
Clearly, \eqref{XformQ} admits a trivial solution if and only if $0= Q(a) = \sqrt{c_1^2 + c_2^2 a^2} = 0, $ or
$c_1 = c_2 = 0$. This proves part (ii) and the theorem. \hfill  $\Box$

\begin{rem}\label{remMp} {\rm Theorem \ref{Xexists} extends (\cite{dou2015}, Thm.~4) from  $\gamma =0$ to $\gamma >0$. A major
shortcoming of Theorem \ref{Xexists} and the above mentioned result in \cite{dou2015} is the presence
of the universal constant $K_p$ 
whose upper bound given in \cite{ose2012} leads to
restrictive conditions on $B_{p,\gamma}$ 
in  \eqref{cQB} for the existence 
of  $L^p$-solution, $p> 2 $. For example, for $p=4$ the above mentioned bound in \cite{ose2012} gives
\begin{equation}\label{K4}
K_4  \mu_4 B_2^2/(1-\gamma)^2 \le (27.083)^4 \mu_4 B_2^2/(1-\gamma)^2 < 1
\end{equation}
requiring $B_2 = \sum_{j=1}^\infty b^2_j $ to be very small. Since statistical inference based of `observable' squares
$r^2_t, 1\le t \le n$
usually requires the existence of 
$\E r^4_t $ and higher moments of $r_t$ (see e.g. \cite{gru2015}), the question arises to derive
less restrictive conditions for the existence of these moments which do not involve the Rosenthal constant
$K_p$. This is achieved in the subsequent Theorem~\ref{Xeven}. Particularly, for $\gamma =0, \lip_Q =1 $  the sufficient condition \eqref{Mp}
of Theorem \ref{Xeven}
for the existence of $\E r^p_t,  p \ge 2 $ even
becomes
\begin{equation}\label{Mp1}
\sum_{j=2}^p {p \choose j}  |\mu_j| \sum_{k=1}^\infty |b_k|^j  \ < \ 1.
\end{equation}
Condition \eqref{Mp1}
coincides with the corresponding condition
in the LARCH case in (\cite{gir2004}, Proposition~3).  
Moreover, \eqref{Mp1} and \eqref{Mp}  apply to more general classes of ARCH models in \eqref{genform0}
and \eqref{genform} to which the specific Volterra series techniques used in  \cite{gir2000}, \cite{gir2004} are not applicable. 
In the particular case $p=4$ condition \eqref{Mp1} becomes
\begin{equation*}
6 B_2 + 4 |\mu_3| \sum_{k=1}^\infty |b_k|^3  + \mu_4 \sum_{k=1}^\infty |b_k|^4  < 1,
\end{equation*}
which seems to be much better than condition \eqref{K4} based on Theorem \ref{Xexists}.
}
\end{rem}

\begin{theo} \label{Xeven} Let $\{\zeta_t\}$ satisfy the conditions of Proposition \ref{Xreq} and
$Q$ satisfy the Lipschitz condition in \eqref{QLip}. \\
Let $p=2,4, \dots $ be even and
\begin{equation}\label{Mp}
\sum_{j=2}^p {p \choose j} |\mu_j| \lip_Q^j  \sum_{k=1}^\infty |b_k|^j   \ < \ (1 - \gamma)^{p/2}.
\end{equation}
Then there exists a unique stationary $L^p$-solution $\{X_t\}$ of \eqref{XformQ}.

\end{theo}

\noi {\it Proof.} For $p=2$, condition \eqref{Mp} agrees with $\lip^2_Q B_{2,\gamma} < 1 $ or condition
 \eqref{cQB} so we shall assume $p \ge 4 $ in the subsequent proof.
In the latter case   \eqref{Mp} implies $\lip^2_Q B_{2,\gamma} < 1 $ and the existence of a stationary $L^2$-solution $\{X_t\}$
of \eqref{XformQ}. It suffices to show that the above $L^2$-solution satisfies $\E X_t^p <  \infty $. \\
Towards this end similarly as in the proof of  Thm \ref{Xexists} (i) consider the solution $\{X^{(n)}_t \}$ with zero initial condition at $t \le -n$
as defined in \eqref{Xndef}. Let $\sigma^{(n)}_t := 0, t < -n$.  Since $\E (X^{(n)}_t - X_t)^2 \to  0 \ (n \to \infty)$,
by Fatou's lemma  it suffices to show that under condition \eqref{Mp}
\begin{equation}\label{Xnbdd}
\E (X^{(n)}_t)^p \  < \  C,
\end{equation}
where the constant $C < \infty $ does not depend on $t, n$. \\
Since $p $ is even for any $t > -n$ we have that
\begin{eqnarray}
\E (X^{(n)}_t)^p
&=&\sum_{s_1,\dots, s_p=-n}^{t-1}  \E \big[ b_{t-s_1} \zeta_{s_1} \sigma^{(n)}_{s_1} \cdots b_{t-s_p} \zeta_{s_p}
\sigma^{(n)}_{s_p} \big]\nn \\
&=&\sum_{j=2}^p {p \choose j} \sum_{s=-n}^{t-1}  b^j_{t-s} \mu_j
\E \Big[  (\sigma^{(n)}_s)^j \Big(\sum_{u=-n}^{s-1} b_{t-u}\zeta_u  \sigma^{(n)}_u \Big)^{p-j} \Big]. \label{Delta1}
\end{eqnarray}
Hence using H\"older's inequality:
$$
|\E \xi^j \eta^{p-j}| \ \le \ c^j \E^{j/p} |\xi/c|^p \E^{(p-j)/p}|\eta|^p \ \le \ c^j \big[ \frac{j}{p c^p}\E |\xi|^p  +
\frac{p-j}{p} \E |\eta|^p\big], \quad  1\le j  \le p, \ c>0
$$
we obtain
\begin{eqnarray}
\E (X^{(n)}_t)^p
&\le&\sum_{j=2}^p {p \choose j} |\mu_j| c_3^j  \sum_{s=-n}^{t-1}  |b^j_{t-s}|  \Big\{ \mbox{$\frac{j}{p c_3^p}$} \E (\sigma^{(n)}_s)^p
+ \mbox{$\frac{p-j}{p}$} \E   \Big(\sum_{u=-n}^{s-1} b_{t-u}\zeta_u  \sigma^{(n)}_u \Big)^p \Big\} \nn \\
&=&\sum_{s=-n}^{t-1}\beta_{1,t-s} \E (\sigma^{(n)}_s/c_3)^p + \sum_{s=-n}^{t-1}\beta_{2,t-s}
\E   \big(X^{(n)}_{t,s}\big)^p, \label{Delta2}
\end{eqnarray}
where $X^{(n)}_{t,s} := \sum_{u=-n}^{s-1} b_{t-u}\zeta_u  \sigma^{(n)}_u,  c_3 > \lip_Q $ and where 
\begin{eqnarray*}
\beta_{1,t-s}&:=&\sum_{j=2}^p  \frac{j}{p} {p \choose j} |b^j_{t-s}|   |\mu_j| c_3^j, \qquad
\beta_{2,t-s}\ := \ \sum_{j=2}^p \frac{p-j}{p}{p \choose j} |b^j_{t-s}|  |\mu_j| c_3^j.
\end{eqnarray*}
The last expectation in \eqref{Delta2}  can be evaluated similarly to \eqref{Delta1}-\eqref{Delta2}:
\begin{eqnarray*}
\E   \big(X^{(n)}_{t,s}\big)^p
&=&\sum_{j=2}^p {p \choose j} \sum_{u=-n}^{s-1}  b^j_{t-u} \mu_j
\E \Big[  (\sigma^{(n)}_u)^j \Big(\sum_{v=-n}^{u-1} b_{t-v}\zeta_v  \sigma^{(n)}_v \Big)^{p-j} \Big] \\
&\le&\sum_{u=-n}^{s-1}\beta_{1,t-u} \E (\sigma^{(n)}_u/c_3)^p
+  \sum_{u=-n}^{s-1}\beta_{2,t-u}
\E  \big(X^{(n)}_{t,u}\big)^p.
\end{eqnarray*}
Proceeding recurrently with the above evaluation results in the inequality:
\begin{eqnarray}
\E (X^{(n)}_t)^p
&\le&\sum_{s=-n}^{t-1} \widetilde \beta_{t-s} \E (\sigma^{(n)}_s/c_3)^p,
 \label{Delta3}
\end{eqnarray}
where
\begin{eqnarray*}
\widetilde \beta_{t-s}&:=&\beta_{1,t-s} \Big(1 + \sum_{k=1}^{t-s-1} \sum_{s < u_k < \dots < u_1 < t} \beta_{2,t-u_1} \cdots \beta_{2,t-u_k}\Big).
\end{eqnarray*}
Let $\beta_i := \sum_{t=1}^\infty \beta_{i,t}, i=1,2,  \, \widetilde \beta := \sum_{t=1}^\infty \widetilde \beta_{t}. $
By assumption \eqref{Mp},
\begin{equation*}
\beta_1 + \beta_2 \ = \  \sum_{j=2}^p {p \choose j} |\mu_j| c_3^j  \sum_{k=1}^\infty |b_k|^j       \   < \  (1-\gamma)^{p/2}
\end{equation*}
whenever $\sigma_3 - \lip_Q >0$ is small enough, 
and therefore
\begin{eqnarray}
\frac{\widetilde \beta}{(1-\gamma)^{p/2}}&\le&\frac{1}{(1-\gamma)^{p/2}} \sum_{t=1}^\infty \beta_{1,t} \big(1 + \sum_{k=1}^{\infty} \beta_2^k) \nn  \\
&=&\frac{1}{(1-\gamma)^{p/2}} \frac{\beta_1}{1- \beta_2} \ < \ 1. \label{unit}
\end{eqnarray}
Next, let us estimate the expectation on the r.h.s. of \eqref{Delta3}
in terms of the expectations on the l.h.s. Using \eqref{Qnel} and Minkowski's inequalities
we obtain
\begin{eqnarray*}
\E^{2/p} (\sigma^{(n)}_s)^p
&\le&\sum_{\ell =0}^{s+n} \gamma^\ell  \E^{2/p} |Q(a+ X^{(n)}_{s-\ell}|^p  \\
&\le&\sum_{\ell =0}^{s+n} \gamma^\ell  \E^{2/p} |c^2_1  + c^2_2(a+ X^{(n)}_{s-\ell})^2|^{p/2} \\
&\le&C +  c^2_3 \sum_{\ell=0}^{n+s}\gamma^\ell \E^{2/p}  (X^{(n)}_{s-\ell})^p,
\label{Delta4}
\end{eqnarray*}
where $c_3 > c_2 > \lip_Q$ and $c_3 - \lip_Q >0$ can be chosen arbitrarily small.
Particularly, for any fixed $ T \in \Z $
\begin{equation*}
\sup_{-n \le s < T} \E^{2/p}(\sigma^{(n)}_s)^p
\ \le\ \frac{c_3^2}{(1-\gamma)}  \sup_{-n \le s < T} \E^{2/p}(X^{(n)}_s)^p  + C. 
\end{equation*}
Substituting the last bound into \eqref{Delta3} we obtain
\begin{eqnarray}
\sup_{-n \le t <T}  \E^{2/p} (X^{(n)}_t)^p
&\le&
\frac{{\widetilde \beta }^{2/p}}{(1-\gamma)}
\sup_{-n \le s < T} \E^{2/p}(X^{(n)}_s)^p
  +  C.
 \label{Delta}
\end{eqnarray}
Relations  \eqref{Delta} and \eqref{unit} imply
\begin{eqnarray*}
\sup_{-n \le t <T}  \E^{2/p} (X^{(n)}_t)^p
\le
 \frac{C}{ 1 - \frac{{\widetilde \beta}^{2/p}}{(1- \gamma)}} \   < \ C
\end{eqnarray*}
proving \eqref{Xnbdd} and the theorem, too. \hfill  $\Box$

\vspace{0.5 cm}

\noi {\bf Example: asymmetric GARCH(1,1). } The asymmetric GARCH(1,1) model of Engle \cite{eng1990} corresponds to
\begin{eqnarray}\label{g1}
\sigma_t^2=c^2+(a+b r_{t-1})^2+\gamma \sigma_{t-1}^2,
\end{eqnarray}
or
\begin{eqnarray}\label{g2}
\sigma_t^2= \theta +\psi r_{t-1}  + a_{11} r^2_{t-1} + \delta \sigma_{t-1}^2
\end{eqnarray}
in the parametrization of (\cite{sen1995}, (5)), with parameters in \eqref{g1}, \eqref{g2} related by
\begin{equation} \label{rel}
\theta=c^2+a^2, \quad \delta=\gamma, \quad \psi=2ab, \quad  a_{11}=b^2.
\end{equation}
Under the conditions that $\{\zeta_t = r_t/\sigma_t\}$ are standardized i.i.d., a stationary asymmetric GARCH(1,1) (or
GQARCH(1,1) in the terminology of \cite{sen1995})  process
$\{r_t\} $ with finite variance and $a\ne 0$ exists if and only if $B_{2,\gamma} = b^2/(1-\gamma) < 1$, or
\begin{equation}\label{bg}
b^2 + \gamma < 1,
\end{equation}
see Thm. \ref{Xexists} (ii). Condition \eqref{bg} agrees with condition $a_{11} + \delta < 1 $ for covariance
stationarity in \cite{sen1995}. Under the assumptions that the distribution of $\zeta_t$ is symmetric and $\mu_4 =  \E \zeta^4_t < \infty $,
\cite{sen1995}
provides a sufficient condition for finiteness
of $\E r^4_t $ together with explicit formula
\begin{eqnarray}\label{Er4}
\E r^4_t&=&\frac{\mu_4 \theta [ \theta (1+a_{11}+\delta) + \psi^2] }{(1-a_{11}^2 \mu_4-2a_{11} \delta -\delta^2) (1-a_{11} - \delta) }.
\end{eqnarray}
The sufficient condition of  \cite{sen1995} for  $\E r^4_t < \infty $
is $\mu_4 a_{11}^2 + 2 a_{11} \delta + \delta^2 < 1$, which translates to
\begin{equation}\label{bg1}
\mu_4 b^4 + 2b^2 \gamma + \gamma^2 < 1
\end{equation}
in terms of the parameters of \eqref{g1}.  Condition \eqref{bg1} seems weaker than 
the sufficient 
condition $\mu_4 b^4 + 6b^2 < (1-\gamma)^2$ of Theorem~\ref{Xeven} for the existence of $L^4$-solution
of \eqref{g1}.  

Following the approach in \cite{dou2015}, below we find explicitly the covariance function $\rho(t) := {\rm cov}(r_0^2, r_t^2)$, including
the expression in \eqref{Er4}, for stationary solution of the  asymmetric GARCH(1,1) in \eqref{g1}.
The approach in \cite{dou2015} is based on derivation and solution of linear  equations for moment functions
$m_2 := \E r_t^2,  \, m_3(t) := \E r_t^2 r_0 $ and $m_4(t) := \E r_t^2 r_0^2 $.
Assume that $\mu_3 = \E \zeta^3_0 = 0$, or
$\E r_t^3=0$. We can write the following moment equations:
\begin{eqnarray} \label{m3}
m_2&=&(c^2+a^2)/(1-b^2-\gamma), \quad m_3(0)=0,\nn \\
m_3(1)&=& \sum_{\ell =0}^{\infty} \gamma^\ell \E (c^2+a^2+2ab r_{-\ell }+b^2 r_{-\ell }^2) r_0=2ab m_2 ,\nn\\
m_3(t)&=& \sum_{\ell =0}^{\infty} \gamma^\ell \E (c^2+a^2+2ab r_{t-\ell -1}+b^2 r_{t-\ell -1}^2) r_0\nn \\
&=& 2abm_2 \gamma^{t-1} + b^2 \sum_{\ell =0}^{t-2} \gamma^\ell m_3(t-\ell -1), \quad t\ge 2.
\end{eqnarray}
From equations above one can show by induction that $m_3(t)=2abm_2 (\gamma+b^2)^{t-1}, t\ge 1$.

Similarly,
\begin{eqnarray*} 
m_4(0)&=&
\mu_4 \E ((c^2+a^2) + 2ab r_0 +b^2 r_0^2 + \gamma \sigma_0^2)^2 \nn \\
&=&\mu_4\Big( (c^2+a^2)^2 + (2ab)^2 m_2+b^4 m_4(0) + 2(c^2+a^2)(b^2+\gamma) m_2 + (2b^2 \gamma +\gamma^2 )m_4(0)/\mu_4\Big),\nn \\
m_4(t)&=&\sum_{\ell =0}^{\infty} \gamma^\ell \E (c^2+a^2+2ab r_{t-\ell -1}+b^2 r_{t-\ell -1}^2) r_0^2\nn \\
&=& \sum_{\ell =0}^{\infty} \gamma^\ell (c^2+a^2) m_2+ b^2 \sum_{\ell =0}^{\infty} \gamma^\ell m_4({|t-\ell -1|}) + 2ab\sum_{\ell =t}^\infty \gamma^\ell m_3({\ell -t+1}), \quad t\ge 1.
\end{eqnarray*}
Using $2ab \sum_{\ell =t}^\infty \gamma^\ell m_3({\ell -t+1})=4 a^2b^2 m_2  \sum_{\ell =t}^\infty \gamma^{\ell} (\gamma+b^2)^{\ell -t}
=4a^2 b^2m_2  \gamma^{t}/(1-\gamma (\gamma+b^2))
$ and $\rho(t)=m_4(t)-m_2^2$ we obtain the system of equations
 \begin{eqnarray} \label{rho4}
 \rho(0) &=&m_4(0)-m_2^2, \nn \\
\rho(t)&=&
b^2 \sum_{\ell =0}^{\infty} \gamma^\ell \rho({|t-\ell -1|}) + 4a^2 b^2m_2  \gamma^{t}/(1-\gamma (\gamma+b^2))\nn \\
&=&
b^2 \sum_{\ell =0}^{t-2} \gamma^\ell \rho({t-\ell -1})
+  C\gamma^{t-1} , \quad t\ge 1,
\end{eqnarray}
where $C:=
b^2 \sum_{\ell = 1}^{\infty} \gamma^\ell \rho({\ell }) +(m_4(0)-m_2^2) b^2+4a^2 b^2 m_2 \gamma/(1-\gamma (\gamma+b^2)) $ is some constant independent of $t$ and
\begin{eqnarray}\label{m4}
m_4(0)&=&
\frac{\mu_4 m_2 }{1-b^4 \mu_4-(2b^2 \gamma +\gamma^2 )}\Big( (c^2+a^2)(1+b^2+\gamma) + (2ab)^2 \Big).
\end{eqnarray}
Note that the expression above coincides with \eqref{Er4} given that the relations in \eqref{rel} hold.

Since the equation in \eqref{rho4} is analogous to \eqref{m3}, the solution to \eqref{rho4} is $\rho(t)=C(\gamma+b^2)^{t-1}, t\ge 1$. In order to find $C$,
we combine $\rho(t)=C(\gamma+b^2)^{t-1}$
and the expression for $C$ to obtain the equation
$
C=
 C b^2 \gamma/(1-\gamma (\gamma+b^2)) +(m_4(0)-m_2^2) b^2+4a^2 b^2 m_2 \gamma/(1-\gamma (\gamma+b^2)).
$
Now $C$ can be expressed as
$$
C=b^2 \frac{(m_4(0)-m_2^2) (1-\gamma(\gamma+b^2))+4a^2 m_2 \gamma}{1-\gamma(\gamma+2b^2)}
$$
together with \eqref{m4} and $\rho(t)=C(\gamma+b^2)^{t-1}, t\ge 1$
giving explicitly the covariances of process $\{r_t^2\}$.
%

\section{Long memory}

The present section studies long memory properties of the generalized quadratic ARCH model in \eqref{genform}
corresponding to $Q(x) = \sqrt{c^2 + x^2}$ of \eqref{Qform1}, viz.,
\begin{equation}  \label{rsqr}
r_t \ = \  \zeta_t \sqrt{\sum_{\ell=0}^\infty \gamma^\ell \big(c^2 + \big(a + \sum_{s<t-\ell} b_{t-\ell-s} r_s\big)^2\big)}, \qquad t \in \Z,
\end{equation}
where $0\le \gamma < 1, a \ne 0, c$ are real parameters,   $\{\zeta_t \}$ are standardized i.i.d. r.v.s, with zero mean and unit variance,
and $b_j, j\ge 1 $ are real numbers satisfying
\begin{equation}\label{Qbc}
b_j \ \sim \ \beta j^{d-1} \ \ (\exists \ 0< d< 1/2, \ \beta>0), 
\end{equation}

The main result of this section is Theorem \ref{long}  which shows that under some additional conditions
the squared process  $\{r^2_t \}$ of \eqref{rsqr}
has similar long memory properties as in case of the LARCH model (see \cite{gir2000}, Thm. 2.2).
Theorem \ref{long} extends the result in (\cite{dou2015}, Thm. 10) to the case $\gamma >0$.
In Theorem \ref{long} and below, $0\le \gamma < 1 $,
$B_2 = \sum_{j=1}^\infty b^2_j$ and $B(\cdot, \cdot)$ is beta function.

\begin{theo} \label{long}
Let $\{r_t\} $ be a stationary $L^2$-solution of \eqref{rsqr}-\eqref{Qbc}.
Assume in addition that
$\mu_4 =\E [\zeta^4_0] < \infty $,
and
$\E [r^4_t ] < \infty $.
Then
\begin{equation} \label{cov2r}
{\rm cov}(r^2_0, r^2_t) \  \sim \  \kappa^2_1 t^{2d-1 }, \qquad  t \to \infty
\end{equation}
where $\kappa^2_1 :=  \big(\frac{2a \beta}{1-\gamma - B_2} \big)^2 B(d, 1-2d) \E r_0^2 $.
Moreover,
\begin{equation}\label{lim2r}
n^{-d-1/2} \sum_{t=1}^{[n\tau]} (r^2_t - \E r^2_t) \ \to_{D[0,1]} \  \kappa_2 W_{d +  (1/2)}(\tau),
 \qquad  n \to \infty,
\end{equation}
where $W_{d +  (1/2)}$ is a fractional Brownian motion with Hurst parameter $H = d + (1/2) \in (1/2,1)$ and
$\kappa^2_2 := \kappa^2_1/ (d(1+2d))$.

\end{theo}

To prove Theorem \ref{long},  we need the following two facts.

\begin{lem} \label{lem1} {\rm (\cite{dou2015}, Lemma 12)} For $\alpha_j \ge 0,$ $j=1,2, \dots, $ denote
\begin{eqnarray*} \label{Ak}
A_k \ := \  \alpha_k + \sum_{0< p <k} \sum_{0< i_1 <  \dots < i_p <  k} \alpha_{i_1} \alpha_{i_2 -i_1}  \cdots \alpha_{i_p - i_{p-1}} \alpha_{k-i_p},
\qquad k =  1,2, \dots.
\end{eqnarray*}
Assume that $\sum_{j=1}^\infty \alpha_j < 1 $ and
\begin{equation*}\label{alpha}
\alpha_j \ \le  \  c \, j^{-\gamma},  \qquad (\exists  \  c >0, \ \gamma >  1).
\end{equation*}
Then there exists  $C >0$ such  that for any $k \ge 1 $
\begin{equation*} \label{Ak1}
A_k  \  \le \   C k^{-\gamma}.
\end{equation*}

\end{lem}

\begin{lem} \label{lem2} Assume that  $0 \le \beta < 1 $ and  $\alpha_j \sim c j^{-\gamma} \  (\exists \, \gamma >0, \, c >0)$. Then
\begin{equation*}
\alpha_{t,\beta} := \sum_{j=0}^{t-1} \beta^j \alpha_{t-j} \ \sim \  \frac{c}{1-\beta}\, t^{-\gamma}, \qquad t \to \infty.
\end{equation*}
\end{lem}

\noi {\it Proof.} It suffices to show that the difference $D_t :=  \alpha_{t,\beta} - \alpha_t/(1-\beta) $ decays faster than $\alpha_t$, in other
words, that
$$
D_t \ = \  \sum_{j=0}^{t-1} \beta^j (\alpha_t -  \alpha_{t-j})  - \sum_{j=t}^\infty \beta^j  \alpha_{t-j} =
o(t^{-\gamma}).
$$
Clearly, $ \sum_{t/2 <j<t} \beta^j (\alpha_t -  \alpha_{t-j}) = O(\beta^{t/2}) = o(t^{-\gamma}), \
\sum_{j=t}^\infty \beta^j  \alpha_{t-j} = O(\beta^{t}) = o(t^{-\gamma})$. Relation
$\sum_{0\le j\le t/2} \beta^j (\alpha_t -  \alpha_{t-j}) = o(t^{-\gamma})$ follows by the dominated
convergence theorem since $\sup_{0 \le j \le t/2} |\alpha_t -  \alpha_{t-j}|t^{\gamma} \le C $ and
$ |\alpha_t -  \alpha_{t-j}|t^{\gamma} \to 0 $ for any fixed $j \ge 0$. \hfill $\Box$

\medskip

\noi {\it Proof of Theorem \ref{long}.} We use the idea of the proof of Thm. 10 in \cite{dou2015}.
Denote
\begin{eqnarray}\label{bgamma}
b_{t,\gamma}&:=&\sum_{j=0}^{t-1} \gamma^j b_{t-j}, \qquad \tilde b^2_{t,\gamma}\ :=\ \sum_{j=0}^{t-1} \gamma^j b^2_{t-j}, \qquad t \ge 1 \\
X_t&:=&\sum_{s<t} b_{t-s} r_s, \qquad  X_{t,\gamma} \ :=\ \sum_{s<t} b_{t-s,\gamma} r_s, \qquad t \in \Z. \nn
\end{eqnarray}
By the definition of $r_t$ in  \eqref{rsqr} we have have the following decomposition (c.f. \cite{dou2015}, (6.66))
\begin{equation}\label{appX}
(r^2_t - \E r^2_t)
- \sum_{s<t} \tilde b^2_{t-s,\gamma} (r^2_s - \E r^2_s) \ = \  2aX_{t,\gamma}  + U_t  + V_{t,\gamma} \ =: \ \xi_t,
\end{equation}
where $X_{t,\gamma}$ is the main term and
the `remainder terms' $U_t$ and $V_{t,\gamma} $ are given by
\begin{eqnarray}
U_t&:=&(\zeta^2_t -\E \zeta_t^2) \sigma^2_t, \qquad V_{t,\gamma}\ :=\  \sum_{\ell =0}^\infty \gamma^\ell V_{t-\ell}, \label{Udef} \\
  \label{Vdef}
V_t&:=&2\sum_{s_2 < s_1  <t} b_{t- s_1} b_{t- s_2} r_{s_1} r_{s_2}.
\end{eqnarray}
Using the identity $V_t = (X^2_t - \E X^2_t) - \sum_{s<t}b^2_{t-2} (r^2_t - \E r^2_t)$ the convergence in $L^2$
of the series on the r.h.s. of \eqref{Vdef} follows as in \cite{dou2015} (6.67).  Hence, the series for $V_{t,\gamma} $
in \eqref{Udef} also converges in $L^2$.

Let us prove that
\begin{equation}\label{gammalim}
{\rm cov}(\xi_0, \xi_t) \sim  4a^2 {\rm cov}(X_{0,\gamma}, X_{t,\gamma}) \sim 4a^2 \lambda^2_1 t^{2d-1}, \ t \to \infty.
\end{equation}
where $\lambda_1^2=\beta^2/(1-\gamma)^2B(d, 1-2d)$.
The second relation in \eqref{gammalim} follows from $b_{t,\gamma} \sim (\beta/(1-\gamma))t^{d-1}, t \to \infty, $ see Lemma \ref{lem1},
and the fact that $X_{t,\gamma} = \sum_{s<t} b_{t-s,\gamma} r_s$ is a moving average in stationary uncorrelated innovations
$\{r_s\}$. Since $\{U_t\}$ is also an uncorrelated sequence, so
${\rm cov}(\xi_0, U_t) = 0 \, (t \ge 1)$, and the first relation in   \eqref{gammalim} is a consequence
of
\begin{eqnarray}
&&\E [U_0 X_{t,\gamma}] + \E [U_0 V_{t,\gamma}]\ =\  o(t^{2d-1}), \label{R1}\\
&&\E [X_{0,\gamma} V_{t,\gamma}] + E [V_{0,\gamma} (X_{t,\gamma} + V_{t,\gamma})]\ =\  o(t^{2d-1}). \label{R2}
\end{eqnarray}
We have $\E [U_0 X_{t,\gamma}] = b_{t,\gamma} \E [U_0 r_0] = O(t^{d-1}) = o(t^{2d-1})$ and
$\E [U_0 V_{t,\gamma}] =  2b_{t,\gamma} D_t  = O(t^{d-1}) = o(t^{2d-1})$, where
$|D_t| := |\E [U_0 r_0 \sum_{s <0} b_{t-s} r_{s}]| \le \E U_0^2 (\E r^4_0)^{1/2} (\E (\sum_{s<0} b_{t-s} r_{s})^4 )^{1/2}  \le C    $
follows from
Rosenthal's inequality in \eqref{rosen} since $\E (\sum_{s<0} b_{t-s} r_{s})^4 \le K_4  \E r^4_0 \big(\sum_{s <0} b^2_{t-s}\big)^2 $
$ \le C$.
This proves  \eqref{R1}. The proof of  \eqref{R2} is analogous to \cite{dou2015} (6.68)-(6.69) and
is omitted.

Next, let us prove \eqref{cov2r}. Recall the definition of $\tilde b^2_{j,\gamma} $ in  \eqref{bgamma}.
From the decomposition  \eqref{appX} we obtain
\begin{equation}\label{r22inv}
r^2_t - \E r^2_t \ = \ \sum_{i=0}^\infty \phi_{i,\gamma} \xi_{t-i},  \qquad t \in \Z,
\end{equation}
where $\phi_{j,\gamma} \ge 0, j\ge 0$ are the coefficients of the power series $\Phi_\gamma(z) := \sum_{j=0}^\infty \phi_{j,\gamma} z^j
=  (1 - \sum_{j=1}^\infty \tilde b^2_{j,\gamma} z^j)^{-1}, \, z \in \C, \, |z| < 1 $ given by $\phi_{0,\gamma} := 1$,
\begin{eqnarray*}\label{phi}
\phi_{j,\gamma}&:=&\tilde b_{j,\gamma}^2+ \sum_{0<k<j} \sum_{0<s_1 < \dots < s_k < j} \tilde b^2_{s_1,\gamma}  \cdots
\tilde b^2_{s_k -s_{k-1},\gamma} \tilde b^2_{j-s_k,\gamma}, \quad j\ge 1.
\end{eqnarray*}
From \eqref{Qbc} and Lemmas \ref{lem1} and \ref{lem2} we infer that
\begin{equation}\label{philim}
\phi_{t,\gamma} =  O(t^{2d-2}), \qquad t \to \infty,
\end{equation}
in particular, $\Phi_\gamma(1) = \sum_{t=0}^\infty \phi_{t,\gamma} = (1-\gamma)/(1- \gamma - B_2) < \infty$ and the r.h.s. of \eqref{r22inv}
is well-defined. Relations \eqref{r22inv} and \eqref{philim}  imply that
\begin{eqnarray}\label{covvv}
{\rm cov}(r^2_t, r^2_0)&=&\sum_{i,j=0}^\infty \phi_i \phi_j {\rm cov}(\xi_{t-i}, \xi_{-j})  \
\sim \ \Phi^2_\gamma(1) {\rm cov}(\xi_{t}, \xi_0), \qquad t \to \infty,
\end{eqnarray}
see \cite{dou2015}, (6.63). Now, \eqref{cov2r} follows from   \eqref{covvv} and  \eqref{gammalim}. The invariance
principle in \eqref{lim2r} follows similarly as in \cite{dou2015}, proof of Thm. 10 from \eqref{r22inv},  \eqref{gammalim} and $n^{-d-1/2} \sum_{t=1}^{[n\tau]} X_{t,\gamma}
\to_{D[0,1]}  \lambda_{2} W_{d +  (1/2)}(\tau) $, $\lambda_{2}^2=\lambda_1^2/d(1+2d)$, the last fact being a consequence of a general result in
\cite{aba2014}. Theorem \ref{long} is proved. \hfill $\Box$

\section{Leverage}

For conditionally heteroscedastic model in \eqref{rsqr} with $\E \zeta_t = \E \zeta^3_t =  0, \E \zeta^2_t = 1 $ consider the leverage function
$h_{t} =  {\rm cov}(\sigma^2_t, r_0) = \E r^2_t r_0, \ t \ge 1. $ Following \cite{gir2004}, and \cite{dou2015},
we say that
$\{r_t\}$ in \eqref{rsqr} {\it has leverage of order $k \ge 1 $}  (denoted by $\{r_t\} \in \ell(k)$) if
\begin{equation*}
h_j < 0, \qquad 1 \le j \le k.
\end{equation*}
The study  in \cite{dou2015} of leverage for model \eqref{rsqr} with $\gamma = 0$, viz.,
\begin{equation*}  \label{rsqr0}
r_t \ = \  \zeta_t \sqrt{c^2 + \big(a + \sum_{s<t} b_{t-s} r_s\big)^2}, \qquad t \in \Z
\end{equation*}
was based on linear equation for leverage function:
\begin{equation*} \label{lev1}
h_t = 2a m_2 b_t + \sum_{0< i < t} b^2_i h_{t-i} + 2 b_t \sum_{i>0} b_{i+t} h_i, \qquad t \ge 1,
\end{equation*}
where $m_2 = \E r^2_0$.
A similar equation \eqref{lev2} for leverage function can be derived for model \eqref{rsqr}  in the general case
$0\le \gamma < 1 $. Namely, using $\E r_s = 0, \,  \E r_s r_0 = m_2 \1(s =0), \E r^2_s r_0 = 0 \, (s \le 0), \E r_0 r_{s_1} r_{s_2} = \1(s_1 =0) h_{-s_2}
\ (s_2 < s_1)$ as in \cite{dou2015} we have that
\begin{eqnarray}
h_t&=&\E r_t^2 r_0 = \sum_{\ell=0}^{t-1} \gamma^{\ell} \E \big[(c^2 + (a+ \sum_{s< t-\ell} b_{t-\ell -s} r_s)^2) r_0\big]\nn \\
&=&\sum_{\ell=0}^{t-1} \gamma^{\ell} \big(2a m_2 b_{t-\ell} +\sum_{s<t-\ell} b^2_{t-\ell-s} \E [r^2_s r_0] \big)
+  2\sum_{\ell =0}^{t-1} \gamma^\ell \sum_{s_2 < s_1 < t-\ell} b_{t-\ell-s_1} b_{t-\ell-s_2}\E [r_{s_1}  r_{s_2} r_0] \nn \\
&=&2a m_2 b_{t, \gamma} +\sum_{0<i<t} h_i \tilde{b}_{t-i,\gamma}^2  + 2\sum_{i>0}  h_i w_{i,t,\gamma},   \label{lev2}
\end{eqnarray}
where $b_{t,\gamma}, \tilde b^2_{t,\gamma} $ are defined in \eqref{bgamma} and $
w_{i,t,\gamma}
:=\sum_{\ell=0}^{t-1} \gamma^{\ell} b_{t-\ell}  b_{i+t-\ell}$.

\begin{prop}\label{lev}
Let $\{r_t\} $ be a stationary $L^2$-solution of \eqref{rsqr} with
$ \E |r_0|^3<\infty,  \, \ |\mu|_3 < \infty$.
Assume in addition that
$ B_{2, \gamma}<1/5, \, \mu_3=\E \zeta_0^3 =0 $.
Then for any fixed $k$ such that $1\le k \le \infty$:

\smallskip

\noi (i) if $a b_1<0$, $ab_j\le 0, j=2, \dots, k$, then $\{r_t\}\in \ell(k)$

\smallskip

\noi (ii) if $a b_1>0$, $ab_j\ge 0, j=2, \dots, k$, then $h_j>0,$ for $ j=1, \dots, k.$
\end{prop}

\noi {\it Proof.} Let us prove that
\begin{equation}\label{hnorm}
\|h\| := \big(\sum_{t=1}^\infty h_t^2\big)^{1/2} \ \le \  \frac{2|a| m_2 B_2^{1/2}}{(1-\gamma)(1-3 B_{2, \gamma})}.
\end{equation}
Let $|b|_{t, \gamma} := \sum_{\ell=0}^{t-1} \gamma^{\ell}  |b_{t-\ell}|$.
By Minkowski's inequality,
\begin{eqnarray}\label{wnorm}
\big(\sum_{i=1}^\infty w^2_{i,t,\gamma}
\big)^{1/2}&\le&\sum_{\ell=0}^{t-1} \gamma^{\ell}  |b_{t-\ell}| \big(\sum_{i=1}^\infty b_{i+t-\ell}^2\big)^{1/2} \ \le \ |b|_{t, \gamma}B_2^{1/2},
\end{eqnarray}
and therefore  $|\sum_{i=1}^\infty  h_i w_{i,t,\gamma} |\le \|h\| B_2^{1/2}|b|_{t, \gamma}$.
Moreover, $\big(\sum_{t=1}^\infty b_{t,\gamma}^2\big)^{1/2} \le B_2^{1/2}/(1-\gamma)$, $\big(\sum_{t=1}^\infty |b|_{t,\gamma}^2\big)^{1/2} $ $\le B_2^{1/2}/(1-\gamma)$ and
$\big(\sum_{t=1}^\infty  (\sum_{0<i<t} h_i \tilde{b}_{t-i,\gamma}^2)^2 \big)^{1/2} \le \|h\| B_{2,\gamma} = \|h\| B_2/(1-\gamma)$.
The above inequalities together  with \eqref{lev2} imply
\begin{eqnarray*}
\|h\|&\le&2|a|m_2B_2^{1/2}/(1-\gamma)+ \|h\| B_2/(1-\gamma) + 2\|h\| B_2^{1/2} B^{1/2}_2/(1-\gamma),
\end{eqnarray*}
proving \eqref{hnorm}.

Using  \eqref{lev2} and \eqref{hnorm},
the statements (i) and (ii) can be proved by induction on $k\ge 1$ similarly to \cite{dou2015}. Since $ w_{i,1,\gamma} =  b_{1}  b_{i+1}$ and $b_{1, \gamma} = b_1$, equation \eqref{lev2} yields
\begin{eqnarray}\label{h1}
h_1&=&2a m_2 b_{1, \gamma}  + 2\sum_{i>0} w_{i,1,\gamma}  h_i =2 b_{1} \big( a m_2 + \sum_{i>0}  h_i   b_{i+1}\big).
\end{eqnarray}
According to \eqref{hnorm}, the last sum in \eqref{h1} does not exceed
$|\sum_{i>0}  h_i   b_{i+1}|\le \|h\| B_2^{1/2}\le  2|a|m_2 B_{2, \gamma}/(1-3 B_{2, \gamma})< |a|m_2$  provided $B_{2, \gamma} <1/5$.
Hence, \eqref{h1} implies  $\text{sgn}(h_1)=\text{sgn}(a b_1)$, or the statements (i) and (ii) for $k=1$.

Let us prove the induction step $k-1 \to k$ in (i). Assume first that $a >0, b_1 <0, b_2 \le 0, \cdots, b_{k-1} \le 0$. Then
$h_1 <0, h_2 < 0, \cdots, h_{k-1} < 0$ by the inductive assumption.  By \eqref{lev2},
\begin{eqnarray*}
h_k&=&2\big(a m_2 b_{k, \gamma}  + \sum_{i>0}  h_i  w_{i, k, \gamma}\big)
+\sum_{0<i<k} \tilde{b}_{i,\gamma}^2 h_{k-i},
\end{eqnarray*}
where $\sum_{0<i<k} \tilde{b}_{i,\gamma}^2 h_{k-i}<0$ and $|\sum_{i>0}  h_i  w_{i, k, \gamma}|\le \|h\|B_2^{1/2} |b|_{k, \gamma}< a m_2 |b|_{k, \gamma}$
according to \eqref{hnorm}, \eqref{wnorm}. Since $b_{k,\gamma} <0 $ and $|b|_{k, \gamma} = |b_{k, \gamma}| $ this implies
$a m_2 b_{k, \gamma}  + \sum_{i>0}  h_i  w_{i, k, \gamma}\le 0$, or $h_k <0$. The remaining cases
in (i)-(ii) follow analogously.  \hfill $\Box$

\section{A simulation study}

As noted in the Introduction,  the (asymmetric) GQARCH model of \eqref{rsqr}  and the LARCH model of \eqref{larch}
have similar long memory and leverage properties
and both can be used for modelling of financial data with the above properties. The main disadvantage of the latter model vs. the
former one seems to be the fact the volatility $\sigma_t$ may take negative values and is not separated from below by positive constant $c>0$
as in the case of \eqref{rsqr}. The standard quasi-maximum likelihood (QMLE)  approach to estimation of
LARCH parameters is inconsistent and other estimation methods were developed
in \cite{ber2009},  \cite{fra2010}, \cite{lev2009}, \cite{tru2014}.   

Consistent QMLE estimation for
5-parametric long memory GQARCH model \eqref{rsqr}  
with $c>0$ and $b_j = \beta j^{d-1} $ was discussed in the recent work \cite{gru2015}.  
The parametric form $b_j = \beta j^{d-1}$ of the moving-average coefficients in \eqref{rsqr} 
is the same as in Beran and Sch\"utzner \cite{ber2009} for the LARCH model.  

It is of interest to compare QMLE estimates and volatility graphs
of the GQARCH and LARCH models based on real data. The comparisons are extended to the classical GARCH(1,1) model
\begin{equation}\label{garch}
r_t = \sigma_t \zeta_t, \quad \sigma_t=\sqrt{\omega  + \alpha r^2_{t-1} + \beta \sigma^2_{t-1}}.
\end{equation}
We consider four data generating processes (DGP):
\begin{eqnarray}\label{DGP}
(L):&&\text{LARCH of \eqref{larch}}, \\
(Q1):&&\text{QARCH of \eqref{rsqr} with $\gamma=0$}, \nn \\
(Q2):&&\text{QARCH of \eqref{rsqr} with $\gamma>0$}, \nn \\
(G):&&\text{GARCH(1,1) of \eqref{garch}}, \nn
\end{eqnarray}
with standard normal innovations and $b_j = \beta j^{d-1}$. The first three models (L), (Q1), (Q2) have long memory and
(G) is short memory.
The parameters
$(a, \beta, d)=(0.0101, -0.1749, $  $0.3520)$ (L), $ (a,c, \beta, d)= (0.0058, -0.0101, 0.2099, 0.4648)$ (Q1),
$(a,c, \beta, d, \gamma)= (0.0020, $  $-0.0049, 0.2394, 0.2393, 0.7735) \  (Q2)$ and
$(\omega, \alpha, \beta) = (0.00001, $
$ 0.1306, 0.8346)\ (G)$
are obtained from real data, consisting of
daily returns of GSPC (SP500) from 2010 01 01 till 2015 01 01 with $n=1257$ observations in total,
by minimizing the corresponding approximate log-likelihood functions.  The details of the estimation procedure 
can be found in \cite{gru2015}.

Fig.~1 presents simulated trajectories of $\sigma_t$ of
four DGP in \eqref{DGP}, corresponding to the same innovation sequence. Observe that the variability of volatility decreases from top
to bottom, (Q2) resembling (G) (GARCH(1,1)) trajectory more closely than
(L) and (Q1).
The graph (Q1) exhibits very sharp peaks and clustering
and a tendency to concentrate near the lower threshold $c$
outside of
high volatility regions. This unrealistic
`threshold effect' is much less pronounced in (Q2) (and also in the other two DGP), due
to presence of the autoregressive parameter $\gamma >0$ which also prevents
sharp changes and excessive variability of volatility series. The graph (G) has different shape and volatility peaks
from the remaining three graphs which is probably due to the short memory of GARCH(1,1).
Fig.~2 illustrates the effect of $\gamma$ on the marginal distribution of (Q2): with $\gamma $ increasing, the
distribution becomes less skewed and  spreads to the right, indicating a less degree of volatility clustering.

\newpage
\vspace{-1 cm}

\begin{figure}[!t]
        \begin{center}
\subfigure{\includegraphics[width=12cm,height=5cm]
                {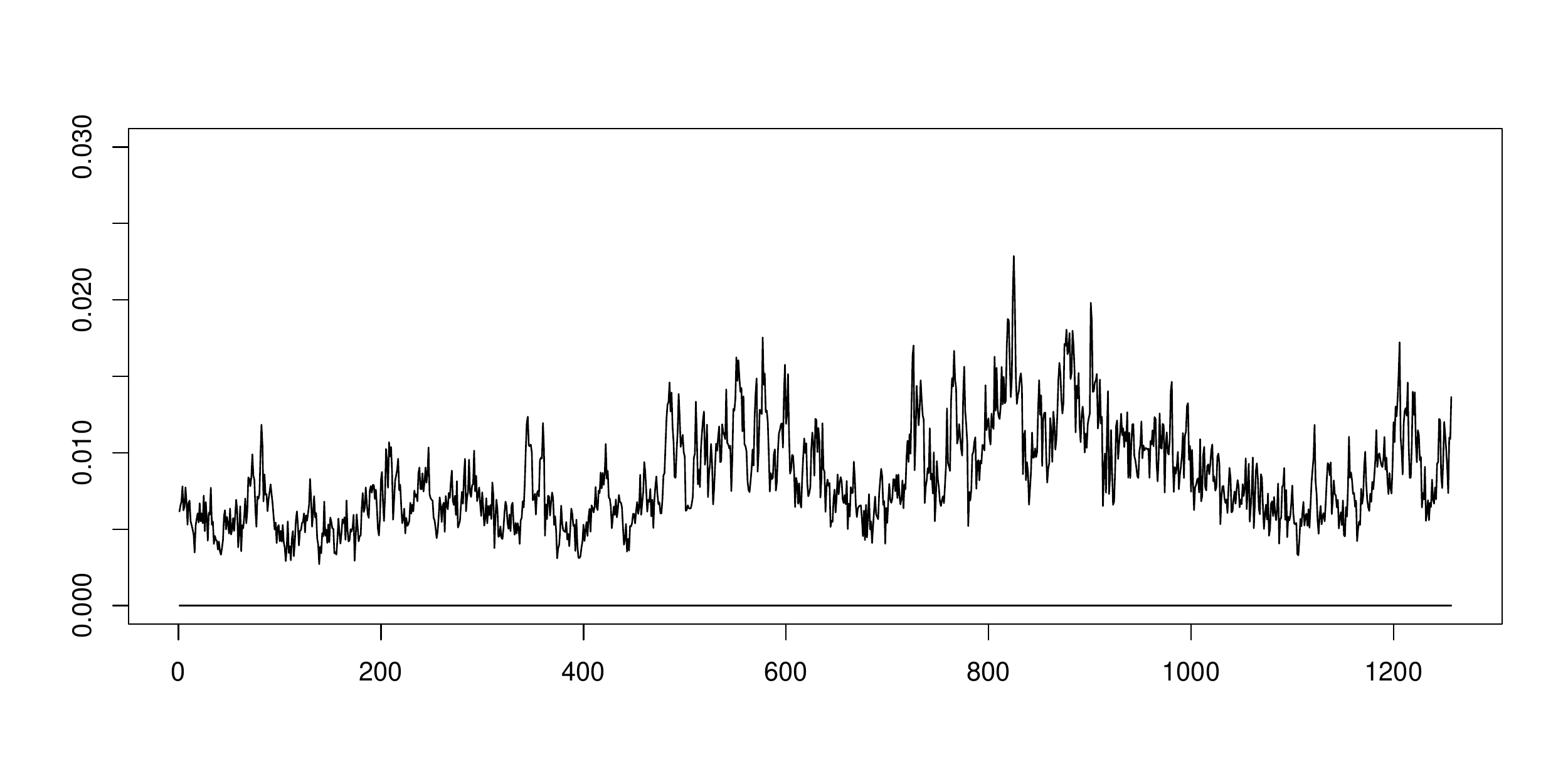}}
                 \vspace{-1.3 cm}

\subfigure{\includegraphics[width=12cm,height=5cm]
                {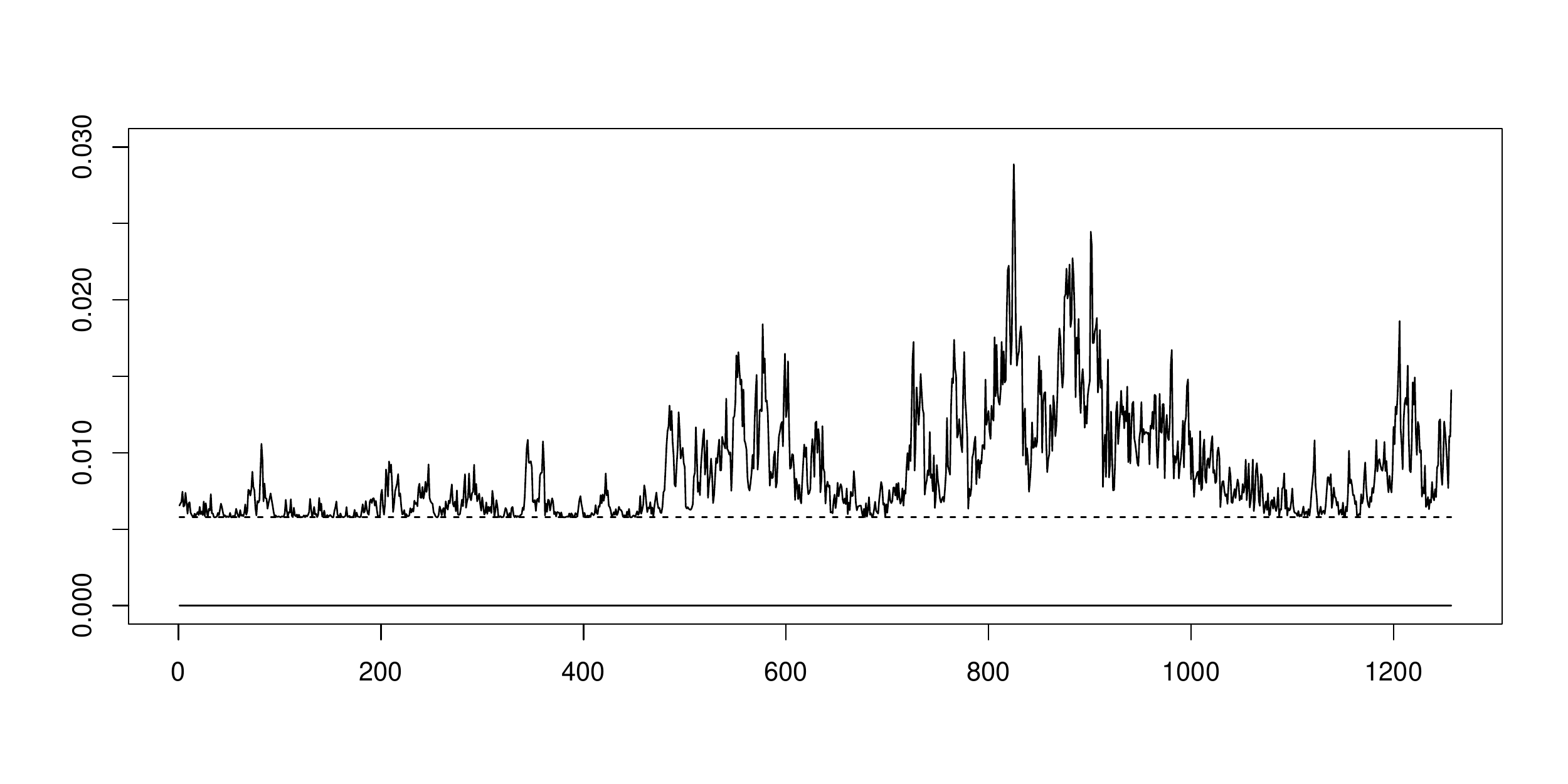}} 	
                 \vspace{-1.3 cm}		%

\subfigure{\includegraphics[width=12cm,height=5cm]
                {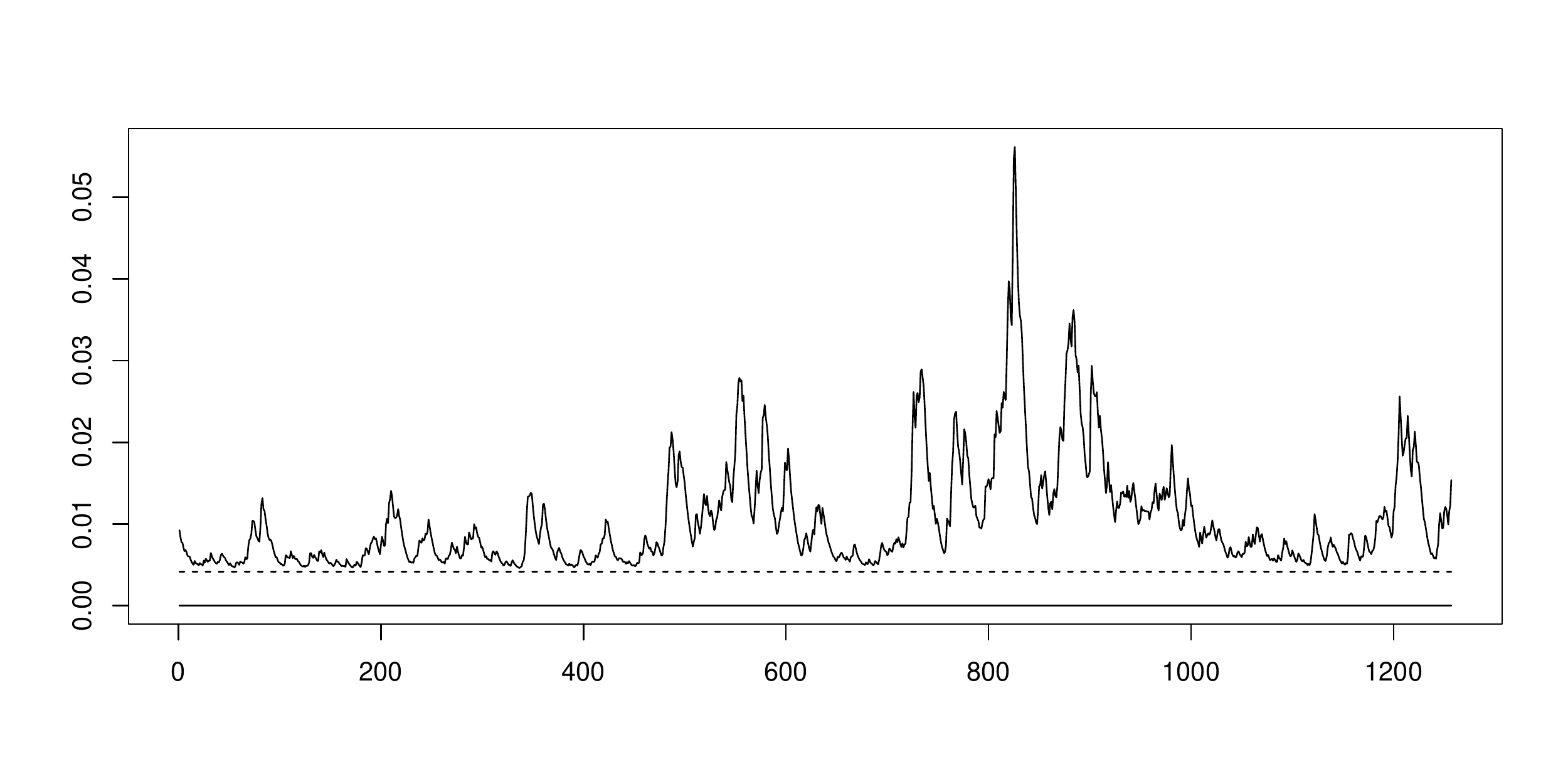}}
                 \vspace{-1.3 cm}

\subfigure{\includegraphics[width=12cm,height=5cm]
                {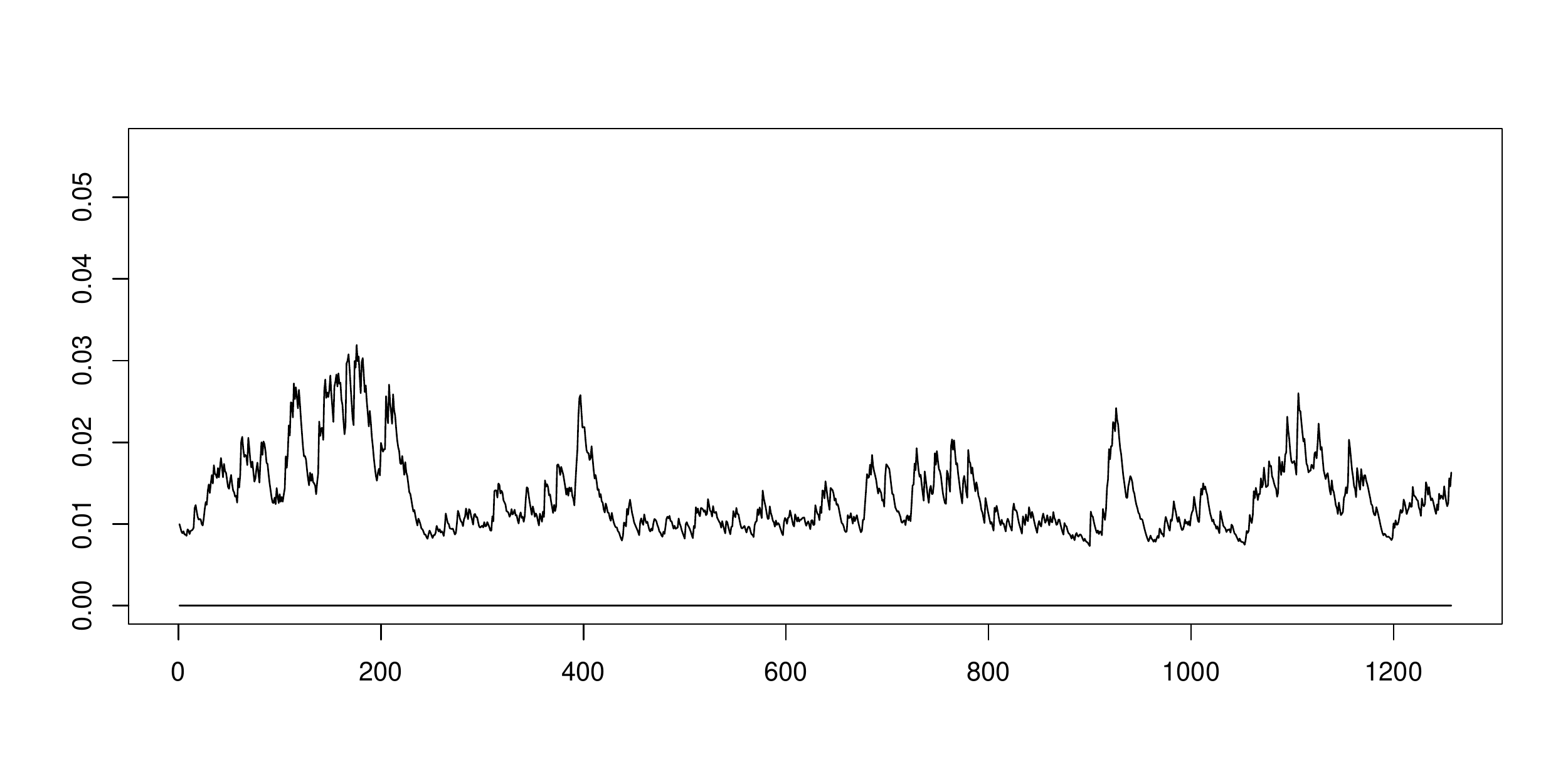}}
                 \vspace{-1.3 cm}

\end{center}
        \caption{\small Trajectory of DGP: From top to bottom: (L), (Q1), (Q2), (G). The dashed line
in  (Q1) and (Q2) indicates the threshold
        $c/\sqrt{1-\gamma}>0$ in \eqref{rsqr}.
        }
\end{figure}

\begin{figure}[!t]
	\centering
	\vspace{-0.7 cm}
	
\includegraphics[height=0.30\textwidth]{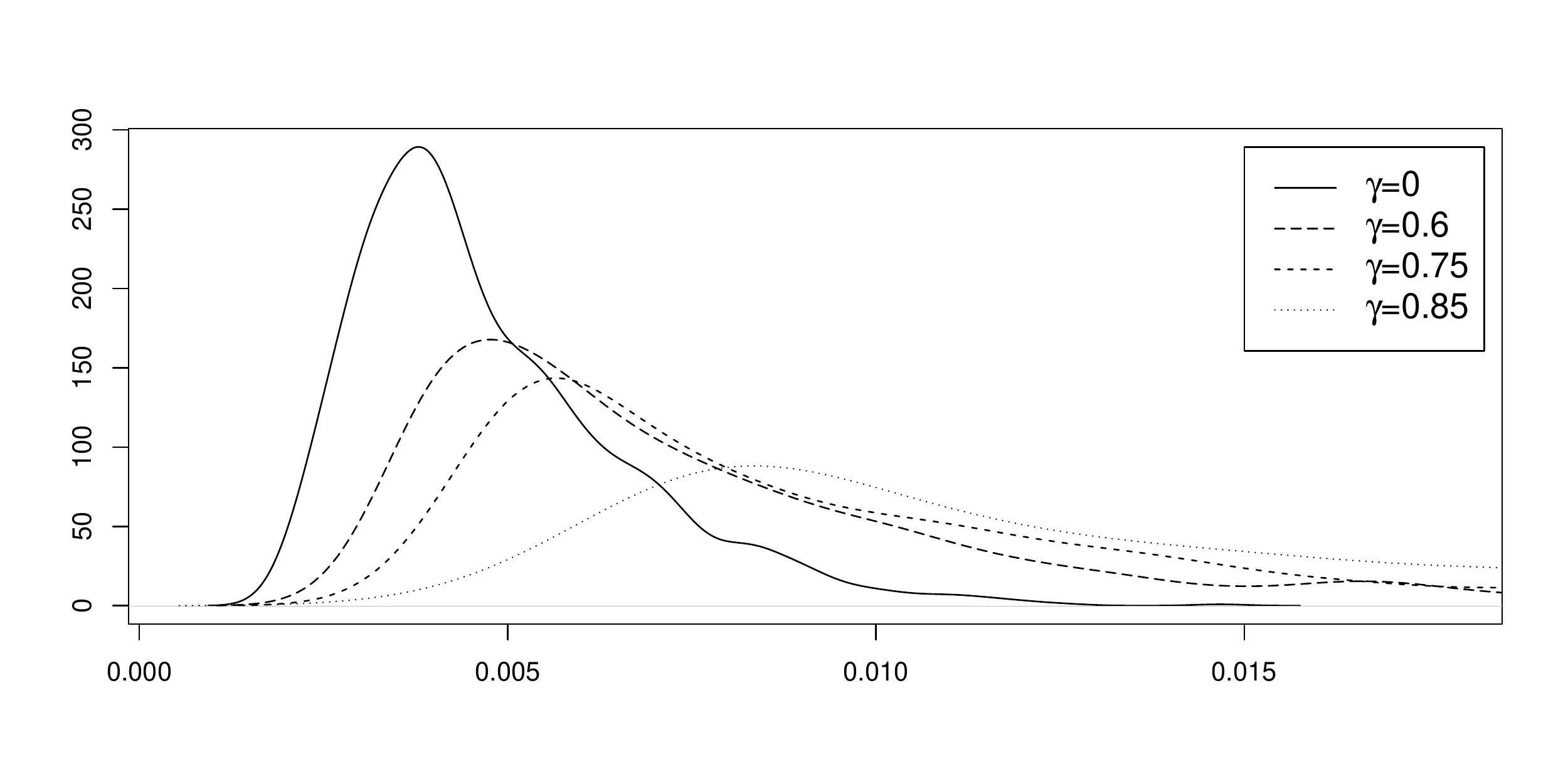}
\vspace{-0.7 cm}

\caption{\small Smoothed histograms of DGP (Q2) for different values of $\gamma$. }
\end{figure}

\end{document}